\documentclass[a4paper,11pt]{article}
\usepackage{amsmath,amscd,amsfonts,amssymb,latexsym}

\newtheorem{theorem}{Theorem}
\newtheorem{proposition}[theorem]{Proposition}
\newtheorem{lemma}[theorem]{Lemma}
\newtheorem{corollary}[theorem]{Corollary}
\newtheorem{remark}[theorem]{Remark}

 \newtheorem{question}[theorem]{Question}

\def\proof{\noindent{\bf Proof.\ }}

\def\qed{~\hbox{$\Box$}\s}

 \def\Q{{\mathbb Q}}
 
\def\Z{{\mathbb Z}}
\def\P{{\mathbb P}}

\def\s{\vskip6pt}

\begin{document}

\title{\bf Purity of boundaries \\of open complex varieties}

\author{ Andrzej Weber\\ \small Department of Mathematics of
Warsaw University\\ \small Banacha 2, 02-097 Warszawa, Poland\\
\small aweber@mimuw.edu.pl}





\date{May 2012}

\maketitle

\begin{abstract} We study the boundary of an open smooth complex algebraic variety $U$.
We ask when the cohomology of the geometric boundary $Z=X\setminus U$ in a smooth compactification $X$ is pure with respect to the mixed Hodge structure.
Knowing the dimension of singularity locus of some singular compactification we give a bound for $k$
above which the cohomology $H^k(Z)$ is pure.
The main ingredient of the proof is purity of the intersection cohomology sheaf.
 \s
\s
    \end{abstract}

\section{Introduction}
Let $U$ be a smooth complex algebraic variety which is not
compact. We study cohomological properties of $U$ which are
invariant with respect to modifications of the interior of $U$. In
other words we investigate cohomolo\-gical properties of the
boundary. The boundary itself can have at least two meanings.
First of all from the topological point of view we may treat an
open smooth variety as the interior of a compact manifold with
boundary. In this case the boundary would mean an odd dimensional
real manifold. We call it the link at infinity. On the other hand
from the geometric point of view we may compactify our variety in
the category of algebraic varieties. The boundary is then a
subvariety of the compactification. In addition we may require
that the compactification is smooth. The condition that the
boundary is a normal crossing divisor is irrelevant for us,
although it is hidden in the construction of the mixed Hodge
structure. The link at the infinity can be identified with the
link of the boundary in the compactification. Regardless from the
differences we show that the topological and geometric boundaries
have a lot in common when the mixed Hodge structure is concerned.

To some extend we try to avoid specific methods of Hodge theory
having in mind possible application (or rather open questions) for
real algebraic geometry, as well as some questions about torsion
for cohomology of complex varieties. The sections
\S\ref{pie}-\S\ref{os} are valid in that generality. Nevertheless
the results of \S\ref{hodz} cannot be generalized and they hold
only for rational cohomology of algebraic varieties. The strong
functoriality of weight filtration implies that lower weight
subspaces of topological and geometric boundary coincide, see
Proposition \ref{teo}. In similar situations this phenomenon was
already described in \cite[Prop. 7.1]{Pa} and \cite[Prop.
5.1]{ABW}.

The proof of the main result of \S\ref{main} uses even stronger
techniques. The purity of intersection cohomology sheaf \cite{BBD}
imposes some conditions on the link of the geometric boundary. We prove
Theorem \ref{teo2} which can be shortened to the following
statement:
\begin{theorem}\label{teo0} Suppose $U$ is a complex smooth algebraic variety.
Assume that $U$ admits a singular compactification $Y$. Suppose
that the singularities of the pair $(Y,Y\setminus U)$ is of dimension $s$. Then for any
smooth compactification $X$ the boundary $X\setminus U$ has pure
cohomology $H^k(X\setminus U)$ for $k\geq \dim(U)+s$.\end{theorem}

The Theorem \ref{teo0} for $U$ admitting one-point compactification
already appeared in \cite[Th.~2.1.11]{dCM}. A vast generalization was
given in \cite{We2}. The present version gives a better bound for
purity, although the situation considered here is less general.

One can treat Theorem \ref{teo0} as a contractibility criterion.
For a subvariety $Z\subset X$ : if $H^k(Z)$ is not pure for $k\geq
\dim(U)+s$, then the pair $(X,Z)$ cannot be contracted to a pair
$(Y,W)$ with singularities of the dimension smaller or equal to
$s$. Although the Theorem \ref{teo0} resembles the Grauert
criterion, it is of different nature. In the Grauert criterion the
intersection form on $Z$ depends on the embedding $Z\subset X$,
whereas here the mixed Hodge structure of $Z$ does not. Of course
it is needles to say that our criterion is not sufficient for
existence of a contraction.

In contrast to the previous paper \cite{We2} we try to present the
subject as elementary as possible. We have avoided to use mixed
Hodge modules \cite{Sa} taking for granted that the cohomology with
coefficients in a complex of sheaves of geometric origin has a
natural mixed Hodge structure. By "geometric origin" we mean
"obtained by the standard sheaf-theoretic operations". We assume
that the varieties are defined over reals or over the complex
numbers and we use classical topology. In our arguments we will
apply resolution of singularities although a part of results
depends only on the formal properties of mixed Hodge modules or
Weil sheaves.

\section{Topological boundary: the link at the infinity}\label{pie}
Let us begin with description of some invariants of open manifolds
which can be defined just using topology and basic properties
resolution of singularities. Working with $\Z/2$-cohomology we can
apply our construction also for real algebraic manifolds. For
complex manifolds we can use any coefficients, not necessarily
$\Q$.

 The first invariant we propose to consider is
the cohomology of the link at infinity:
 $$H^*(L_\infty U ):=\lim_{\overrightarrow{K\subset  U}}H^*( U\setminus K)\,,$$
where $K$ runs through compact sets contained in $U$. The group
$H^*(L_\infty U )$ is exactly the cohomology of the link $L_Z$,
the link of the boundary set $Z=X\setminus U$, where $X$ is a compactification of $U$. (For various
approaches to the link of a subvariety, see \cite{DS}.) This
cohomology group is of finite dimension. It can be expressed in
terms of sheaf operations on $X$:
 $$H^*(L_\infty U
)=H^*(L_Z )=H^*(Z;i^*Rj_*\Q_U)\,,$$
 where $j:U\hookrightarrow X$ and $i:X\setminus U\hookrightarrow X$ are the
inclusions.

\section{Geometric boundary: image of boundary cycles}

Another invariant considered by us is the image of boundary cycles
$$IB^*(U)=
{\rm im}(H^*(Z)\to H^{*+1}_c(U))={\rm ker}(H^{*+1}_c(U)\to
H^{*+1}(X))\,,$$
 where $X$ is a {\bf smooth} compactification of
$U$ and $Z=X\setminus U$. The maps come from the long exact
sequence of the pair $(X,Z)$. \s

To show the independence of $X$ we start with a purely topological lemma.

\begin{lemma} \label{lem} Suppose we have a map of real smooth oriented closed
manifolds
 $$f:X_1\to X_2$$
 which is isomorphism of some open subsets
$$f_{|U_1}:U_1=f^{-1}(U_2)\stackrel{\simeq}{\to}U_2\,.$$
 Denote by
$$IB_i^*={\rm ker}(H^{*+1}_c(U_i)\to
H^{*+1}(X_i))$$ the kernels of the natural maps for $i=1,2$. Then $f$
induces an isomorphism
$$f^*:IB_2^*\to IB_1^*\,.$$
\end{lemma}

\proof The map $f^*:H^*(X_2)\to H^*(X_1)$ is injective since it is
a map of degree one of compact manifolds. The map $f$ induces the
transformation
$$\begin{matrix}
 IB_2^k &\hookrightarrow &H^{k+1}_c(U_2)&\to& H^{k+1}(X_2)\cr
 \downarrow&&^\simeq\downarrow&&^{mono}\downarrow\cr
 IB_1^k &\hookrightarrow &H^{k+1}_c(U_1)&\to& H^{k+1}(X_1)\cr
\end{matrix}$$
It follows that $IB_2^k\to IB_1^k$ is an isomorphism.\hfill\qed

To prove the independence of $IB^*(U)$ on the compactification it
remains to say that any two smooth compactifications are dominated
by a third one.

\section{Basic exact sequences}\label{os}

We will need three exact sequences to relate the described
invariants. These exact sequences may be constructed
topologically, but it is important to know that they come from
distinguished triangles in the derived category of sheaves. It
will follow, that for complex varieties the maps of the described
exact sequences preserve the mixed Hodge structure. \s

We start with the sequence relating the cohomology of $U$ and the
cohomology of its link at the infinity. Let $X$ be any
compactification and $Z=X\setminus U$.
 It is possible to find a neighbourhood $N$ of $Z$
which retracts to $Z$ and the boundary  $\partial N$ is
homeomorphic to the link of $Z$.
 Considering the pair $(X\setminus N,\partial N)$ we arrive to the
long exact sequence
\begin{equation}\label{piecio}\to H^k(U)\to H^k(L_Z)\stackrel\delta\to H^{k+1}_c(U)\to H^{k+1}(U)\to\,.\end{equation}
This exact sequence may be in fact  obtained from the fundamental
distinguished triangle (in the category of mixed Hodge modules on
$X$)
\begin{equation}\label{trojkat}\begin{matrix} i_*i^!G&\longrightarrow&\phantom{gg} G\cr
      _{[+1]}\hfill    \nwarrow&&\swarrow\hfill&\cr
          & Rj_*j^*G \cr\end{matrix}\end{equation}
where $G=j_!\Q_U$. By duality we obtain the triangle
$$\begin{matrix} i_!i^*Rj_*\Q_U&\longleftarrow&\phantom{gg} Rj_*\Q_U\cr
      _{[+1]}\hfill    \searrow&&\nearrow\hfill&\cr
          & j_!\Q_U\cr\end{matrix}$$
since $j_!j^!Rj_*\Q_U\simeq j_!\Q_U$. Applying the cohomology we
obtain the sequence (\ref{piecio}). \s

We also need an exact sequence relating $H^k(Z)$ and $H^k(L_Z)$.
Topolo\-gi\-cally we have  a retraction  $N\to Z$. The exact
sequence for the manifold with boundary $(N,\partial N)$
$$
 \to H^{k}(N)\to H^{k}(\partial N)\to H^{k+1}(N,\partial N)\to H^{k+1}(N)\to $$
 becomes
 \begin{equation}\label{egse}
 \to H^{k}(Z) \to  H^{k}(L_Z) \to  H^{k+1}(X,U) \to  H^{k+1}(Z) \to  
 \,.\end{equation}
 The sheaf
theoretic definition is given below. Let us restrict the triangle
(\ref{trojkat}) with $G=\Q_X$ to $Z$. We have
$i^*i_*i^!\Q_X=i^!\Q_X$ and we obtain the triangle
$$\begin{matrix} i^!\Q_X&\longrightarrow&\Q_Z\cr
      _{[+1]}\hfill    \nwarrow&&\swarrow\hfill&\cr
          & i^*Rj_*\Q_U.\cr\end{matrix}$$
The associated sequence of cohomology is just (\ref{egse}). It
 plays the  fundamental role in our further consideration.

\s Of course the third exact sequence used by us is the sequence
of the pair $(X,Z)$
 \begin{equation}\label{trzeci}\to H^k(X)\to H^k(Z)\to H^{k+1}_c(U)\to
H^{k+1}(X)\to\,.\end{equation}
 To relate the groups $H^*(L_\infty(U))=H^*(L_Z)$ and
 $IB^*(U)={\rm im}(H^k(Z)  \to  H^{k+1}_c(U))$ we
  apply the map of exact sequences (\ref{piecio}) and
 (\ref{trzeci}) induced by the inclusion
 $$(X\setminus int(N), \partial N)\subset (X,N)\,.$$
 We obtain the commutative diagram
 $$\begin{matrix} H^k(Z) & \to & H^{k+1}_c(U)  \cr
 \downarrow && \parallel \cr
 H^k(L_Z) & \to & H^{k+1}_c(U) &. \cr\end{matrix}$$
 We see that
 $$IB^k(U)= {\rm im}(H^k(Z) \to
 H^{k+1}_c(U))\subset {\rm im}(H^k(L_\infty(U)) \to
 H^{k+1}_c(U))\,.$$
In general the inclusion is proper.

\section{Mixed Hodge structure}\label{hodz}
From now on we consider only complex algebraic varieties and
rational cohomology.

The considered invariants $H^k(L_\infty U)$ and $IB^k(U)$
 are equipped with a mixed Hodge structures. The first one is given by the sheaf-theoretic description:
  $$H^*(L_\infty U
)=H^*(Z;i^*Rj_*\Q_U)\,.$$
 The second one, $IB^*(U)$, has a structure induced from $H^*(Z)$. In the situation of Lemma \ref{lem}
 the map $IB_2^k\to IB_1^k$ preserves quotient mixed Hodge structures and
since it is an isomorphism of vector spaces it must be also an
isomorphism of all weight subspaces. In fact by the definition of
the mixed Hodge structure we have
$$IB^k(U)=W_kH^{k+1}_c(U)\,.$$

For us the most interesting part is the weight subspace $W_{k-1}$.
Using basic properties of the mixed Hodge structure we will give
three description of that weight space.

\begin{proposition}\label{teo} Let $X$ be a smooth compactification of $U$
and $Z=X\setminus U$. Then the following groups are isomorphic:
\begin{enumerate}
\item $ W_{k-1}H^k(L_\infty U)\,,$
\item $W_{k-1}H^{k+1}_c(U)\,,$
\item $W_{k-1}H^k(Z)\,.$
\end{enumerate}
\end{proposition}
Note that in the statement of the theorem we do not assume that
$Z$ is a smooth divisor with a normal crossing.
 As a corollary from Proposition
\ref{teo} we have

\begin{corollary} Let $X$ be a smooth compactification of $U$ and $Z=X\setminus U$.
 The cohomology $H^k(Z)$ is pure of weight $k$ if and only if
$H^{k}(L_Z)$ is of weight $\geq k$.
\end{corollary}

Also we note (compare \cite[Prop. 7.1]{Pa}):
\begin{corollary}\label{cor1}
 The impure part of cohomology of the boundary set
$W_{k-1}H^k(Z)$ does not depend on the smooth compactification.
\end{corollary}

\begin{remark}\rm Note that the group  $W_{k -1}H^k (Z)$ is a
topological invariant of $Z$, since by \cite{We} it is the kernel
of the canonical map to the intersection cohomology $H^k(Z)\to
IH^k(Z)$. Also by the construction of the mixed Hodge structure
we have
$$W_{k-1}H^k(Z)=ker(g^*:H^k(Z)\to H^k(\widetilde Z))\,,$$
where $g:\widetilde Z\to Z$ is any dominating proper map from a
smooth variety, possibly of bigger dimension. \s
\end{remark}

The entire cohomology of the boundary of a smooth compactification
is not an invariant of $U$. Of course when we blow up something at
the boundary then the cohomology is modified, nevertheless the lower
parts of weight filtration remains unchanged.

\begin{remark}\rm With help of the Decomposition Theorem of \cite{BBD} we have better
insight to what happens with the cohomology of the boundary. Let
$f$ be a map of pairs $(X_1,Z_1)\to (X_2,Z_2)$ which is an
isomorphism outside $Z_1$. The push-forward of the constant sheaf
on $X_1$ decomposes:
$$Rf_*\Q_{X_1}\simeq\Q_{X_2}\oplus \bigoplus_\alpha IC(V_\alpha;L_\alpha)\,.$$
The supports of the intersection sheaves $IC(V_\alpha;L_\alpha)$
are contained in $Z_2$, therefore
$$H^*(Z_1)=H^*(Z_2;(Rf_*\Q_{X_1})_{|Z_2})
\simeq H^*(Z_2)\oplus \bigoplus_\alpha
IH^*(V_\alpha;L_\alpha)\,.$$ Again we see that the difference
between $H^*(Z_1)$ and $H^*(Z_2)$ is pure since $\bigoplus_\alpha
IH^*(V_\alpha;L_\alpha)$ is a summand of $H^*(X_1)$ .\end{remark}

\begin{remark}\rm Using another powerful tool, namely the Weak Factorization
Theorem \cite{AKW}, we can trace how the cohomology of the
boundary may change. Each time when we blow up a smooth center $S$
contained in the boundary the pure summand ~${\rm
coker}\!\left(H^*(S)\to H^*(\P N_{S/X})\right)$ contributes to the
cohomology of the blown up boundary. Here $H^*(\P N_{S/X})$ is the
projectivization of the normal bundle of $S$ in $X$.
\end{remark}

 The proof of Proposition \ref{teo} is divided into Lemmas \ref{p1} and \ref{p2}.

\begin{lemma} \label{p1} We have $$W_{k-1}H^k(Z)\simeq W_{k-1}H^{k+1}_c(U)\,.$$
\end{lemma}

\proof We recall that $H^k(X)$ is of weight $k$ and  $H^{k+1}(X)$ is of weight $k+1$.  Therefore the long exact sequence
$$\to H^k(X)\to H^k(Z)\stackrel\delta\to H^{k+1}_c(U)\to H^{k+1}(X)\to$$
induces an isomorphism of graded pieces for $\ell<k$
$$Gr_\ell^WH^k(Z)\simeq Gr_\ell^W H^{k+1}_c(U)\,.$$
It follows that the boundary map $W_{k-1}H^k(Z)\to
W_{k-1}H^{k+1}_c(U)$ is an isomorphism.\hfill\qed

\begin{lemma} \label{p2} We have $$W_{k-1}H^k(L_Z)\simeq W_{k-1}H^{k+1}_c(U)\,.$$
\end{lemma}

\proof We consider the long exact sequence (\ref{piecio}). Since
$U$ is smooth $W_{k-1}H^k(U)=0$. Therefore for $\ell<k$
$$Gr_\ell^WH^k(L_Z)\simeq Gr_\ell^W H^{k+1}_c(U)\,.$$
Again the boundary map $W_{k-1}H^k(L_Z)\to W_{k-1}H^{k+1}_c(U)$ is
an isomorphism. \hfill\qed

\section{Singular versus smooth compactifications}\label{main}

Let $W\subset Y$ be a pair of varieties. Assume that $Y\setminus
W$ is smooth. By the singularity of the pair we mean the set of
points at which $W$ in $Y$ analytically does not look like a
submanifold (of any dimension) in a manifold. The singularity set
consists of points at which  $W$ or $Y$ is singular.
 Below we give the exact statement of our main result.

\begin{theorem} \label{teo2} Let  $U$ be a smooth variety. Suppose that $U$  admits a compactification $Y$ and let
$W=Y\setminus U$ be the boundary set. Denote by $s$ the dimension
of the singularities of the pair $(Y,W)$. Let $X$ be a smooth
compactification of $U$ and $Z=X\setminus U$.
 \noindent For $k\geq \dim(U)+\dim(W)$ we have:

 i) the cohomology of
the link $H^{k}(L_Z)$ is of weight $\geq k+1$,

ii) the restriction map $H^k(Z)\to H^k(L_Z)$ vanishes.

\s \noindent For $k\geq \dim(U)+s$ we have:

 iii) the cohomology
of the boundary $H^k(Z)$ is pure of weight $k$, that is
$$W_{k-1}H^k(Z)=0\,$$

iv) the cohomology of the link $H^{k}(L_Z)$ is of weight $\geq k$.
\end{theorem}

Note that by  Proposition \ref{teo} the  claim iv) does not depend
on the choice of the smooth compactification $X$.

\s Let $n=\dim(U)$. By Poincar\'e duality we have
$$H^k(Z)^*=H^{2n-k}(X,U)(n)\,,$$
$$H^k(L_Z)^*=H^{2n-1-k}(L_Z)(n)\,,$$
where $(n)$ denotes the Tate twist shifting the weights by $2n$.
The dual version of the Theorem \ref{teo2} is the following:
\begin{theorem}\label{teodu} With the assumption of Theorem
\ref{teo2}:

\s \noindent For $k\leq \dim(U)-\dim(W)$  we have

i') the cohomology of the link $H^{k-1}(L_Z)$ is of weight $\leq
k-1$,

ii') the boundary map $H^{k-1}(L_Z)\to H^k(X,U)$ vanishes.

\s\noindent For $k\leq  \dim(U)-s$  we have

iii') the cohomology $H^k(X,U)$ is pure of weight $k$, that is
$$W_{k}H^k(X,U)=H^k(X,U)\,$$

iv') the cohomology of the link $H^{k-1}(L_Z)$ is of weight $\leq
k$.
\end{theorem}

To distinguish two copies of $U$ in $X$ and in $Y$ we will use the
letter $V$ for the copy of $U$ in $Y$. The identification map
$U\to V$ is denoted by $f$:
$$\begin{matrix} Z=X\setminus U &\subset & X=\overline U & \supset& U\cr\cr &&&^\simeq&\big\downarrow&\hskip-5pt
^f\cr\cr W=Y\setminus V &\subset & Y =\overline V& \supset &V\cr
\end{matrix}$$

\begin{remark} \label{comp} \rm  In our setup, we
can apply completion and resolution of singularities. Therefore
$X$ can be replaced by a dominating smooth variety for which the
map $f$ extends to the boundary.
\end{remark}

Some information about the weights of cohomology of the link and
the boundary  can be deduced when we have a proper map $f:U\to V$
and a compactification of $V$. A statement which generalizes i)
and ii) in terms of a defect of semismallness \cite{dCM} is
formulated in \cite{We2}.
 The direct generalization of iii) and iv) would involve precise
 information about the singularities of the perverse cohomology sheaves
 $^p{\cal H}^kRf_*\Q_U$.

\s

The Theorem \ref{teo2} can be localized around a topological
component of $X\setminus U$. Precisely, consider the set of ends,
i.e. $U_\infty=\pi_0( X\setminus  U)$. This set does not depend on
the choice of $X$ provided that $X$ is normal.
 A  map of algebraic varieties  which is proper induces a map of their
ends. To deduce purity of the cohomology of a part of the boundary
of $U$ it is enough to have information about a singular
completion of  the corresponding end.

\section{Proofs}

 Before the proof of Theorem \ref{teo2} let us recall the key property of the link of a subvariety

\begin{theorem}[\cite{DS}] \label{desa} Let $Y$ be a variety and let $W$ be a compact subvariety. Let us assume that
$Y\setminus W$ is smooth. Then $H^k(L_W)$ is of weight $\leq k$
for the degrees $k<\dim(Y)-\dim(W)$.\end{theorem}

Theorem \ref{desa} immediately follows from the purity of the
intersection sheaf \cite{G,BBD} since the stalk cohomology ${\cal
H}^k(IC_Y)$ is isomorphic to ${\cal H}^k(Rj_*\Q_V)$ for
$k<\dim(Y)-\dim(W)$ and $H^*(L_W)=H^*(W;(Rj_*\Q_V)_{|W})$.

\begin{remark}\rm In \cite[\S6]{ABW} the Decomposition Theorem of \cite{BBD} was used  to give estimates for the dimension of intersection cohomology of the link by means of resolution. But it seems that the purity of the intersection sheaf was not used directly. \end{remark}

\s {\bf Proof of (\ref{teodu}.i'-ii')}.

By Remark \ref{comp} we assume that the map $f$ extends to $X$.
The extended map  (denoted by the same letter) induces a map of
sheaves $i'^*Rj'_*\Q_V\to Rf_*i^*Rj_*\Q_U$. Therefore the mixed
Hodge structures of the isomorphic groups $H^*(L_W)$ and
$H^*(L_Z)$ coincide. By Theorem \ref{desa} and the assumption on the
dimension of $W$ the cohomology $H^{k-1}(L_W)$ is of weight $\leq
k-1$.
 The claim \ref{teodu}.ii' follows from the long exact sequence
(\ref{egse}): the boundary map
$$H^{k-1}(L_Z)\to H^k(X,U)$$ vanishes because the first term is of
weight $\leq k-1$ and the second term is of weight $\geq k$.\s

{\bf Proof of (\ref{teodu}.i-ii)} follows by duality.\s

{\bf Proof of (\ref{teodu}.iii-iv)}   If $W=Sing(Y)$ then
$s=\dim(W)$ and the statement  i) is even stronger then required.
The Proposition \ref{teo}  implies  ii).\s

 Suppose now $Sing(Y)\subsetneq W$. We may assume that $f$
extends to a map $X\to Y$ and also we may assume that the map $f$
is a resolution of singularities of the pair $(Y,W)$. Let
$\widetilde W\subset X$ be the proper transform of $W$. Denote by
$E\subset X$ the exceptional set of $f$ and let $F=E\cap
\widetilde W$. Consider the Mayer-Vietoris exact sequence for
$Z=E\cup \widetilde W$:
$$\to H^{k-1}(E)\oplus H^{k-1}(\widetilde W)\stackrel\alpha\to H^{k-1}(F)\stackrel\delta\to H^k(Z)\to H^k(E)\oplus
H^k(\widetilde W)\to.$$ By
 (\ref{teo2}.i) applied to the map $(X,E)\to (Y,f(E))$ the
 cohomology of the link $H^k(L_E)$ is of weight $\geq k+1$ for $k\geq
 \dim(X)+s$. Hence by Proposition \ref{teo}
the cohomology $H^k(E)$ is pure for $k\geq \dim(X)+s$.  Of course
$H^k(\widetilde W)$ is pure since we assume that $\widetilde W$ is
smooth. To prove the purity of $H^k(Z)$ it remains to show that
the map $\delta$ of the Mayer-Vietoris sequence is trivial.
\s
By (\ref{teo2}.ii) applied to $F\subset\widetilde W$ the map
$$H^{k-1}(F)\to H^{k-1}(L_F)$$ vanishes for $k-1\geq\dim(W)+s$. By
the exact sequence (\ref{egse}) for that pair the restriction map
 $$H^{k-1}(\widetilde W,\widetilde W\setminus F)\to H^{k-1}(F)$$
 is surjective.
  The above map factors through $H^{k-1}(\widetilde W)$, therefore the  map
   $$H^{k-1}(\widetilde W)\to H^{k-1}(F)$$
   is surjective.    It follows that the restriction map
 $\alpha$ is surjective and the boundary map $\delta$ is
  trivial for $k\geq\dim(X)+s\geq \dim(W)+s+1$. This completes the proof.\s

{\bf Proof of (\ref{teodu}.iii'-iv')} follows by duality.
 \hfill\qed

\begin{remark}\rm If the singularity set is empty then $s=-\infty$ by convention. The claims (\ref{teo2}.i-ii) hold for all degrees by trivial reasons.\end{remark}

The special case when $W$ is a point (an isolated singularity
resolution) was studied from the very beginning of the theory. In
that case both maps $H^{n-1}(L_Z)\to H^n(X,U)$ and $H^n(X)\to
H^n(L_Z)$ are trivial. The map $H^n(X,U) \to H^n(Z)$ is an
isomorphism. After the identification $H^n(X,U)=H^n(Z)^*$ we
obtain a nondegenerate intersection form which was studied for
example in \cite{GM}.

\section{Questions about real algebraic varieties}
The Hodge theory for real algebraic varieties and $\Z/2$
coefficients is not available. The approach of \cite{MCPa1, MCPa2}
does not lead to a strongly functorial weight filtration.
Nevertheless one defines impure cohomology of a singular compact
variety $X$: it is the kernel of $H^*(X)\to H^*(\tilde X)$, where
$\tilde X$ is any resolution. We say that the cohomology of a real
variety is pure if the kernel  $H^*(X)\to H^*(\tilde X)$ is
trivial. The definition does not depend on $\tilde X$. One can ask
the question about the generalization of the Theorem \ref{teo2}:

\begin{question} With the assumption of Theorem \ref{teo}
 for real algebraic varieties: What properties of $(Y,W)$ would guarantee
purity of $H^*(Z)$ in some range of degrees?\end{question}

The dimension of the singularity set is far to weak invariant. It
is well known that any real algebraic set can be contracted to a
point.

 \end{document}